\documentclass[12pt,a4paper]{article}
\usepackage{amssymb,amsmath,amsthm,color, enumerate, soul}

\newcommand\cC{{\mathcal C}}

\newcommand\cF{{\mathcal F}}

\newcommand\cH{{\mathcal H}}

\newcommand\cK{{\mathcal K}}

\newcommand\cM{{\mathcal M}}

\newcommand\cP{{\mathcal P}}

\newcommand\cS{{\mathcal S}}

\theoremstyle{plain}
\newtheorem{theorem}{Theorem}[section]
\newtheorem{lemma}[theorem]{Lemma}
\newtheorem{claim}[theorem]{Claim}
\newtheorem{observation}[theorem]{Observation}
\newtheorem{corollary}[theorem]{Corollary}

\theoremstyle{definition}

\newtheorem{fact}[theorem]{Fact}

\newtheorem{remark}[theorem]{Remark}

\newcommand\sref[1]{Section~\ref{sec:#1}}



\title{Avoider-Enforcer star games}
\author{Andrzej Grzesik \thanks{Theoretical Computer Science Department,
Faculty of Mathematics and Computer Science, Jagiellonian
University, ul. Prof. St. Lojasiewicza 6, 30-348 Krakow, Poland.
Email:andrzej.grzesik@uj.edu.pl }  \and Mirjana Mikala\v cki
\thanks{Department of Mathematics and Informatics, Faculty of Sciences, University of
Novi Sad, Serbia. Research partly supported by Ministry of Education
and Science, Republic of Serbia, and Provincial Secretariat for
Science, Province of Vojvodina. Email:
mirjana.mikalacki@dmi.uns.ac.rs. }  \and Zolt\'an  L\'or\'ant Nagy
\thanks{Alfr\'ed R\'enyi Institute of Mathematics, P.O.B. 127,
Budapest H-1364, Hungary. Supported by Hungarian National Scientific
Research Funds (OTKA) grant 81310. Email:
nagy.zoltan.lorant@renyi.mta.hu.}  \and Alon Naor \thanks{School of
Mathematical Sciences, Raymond and Beverly Sackler Faculty of Exact
Sciences, Tel Aviv University, Tel Aviv, 69978, Israel. Email:
alonnaor@post.tau.ac.il}  \and Bal\'azs Patk\'os \thanks{MTA--ELTE Geometric and 
Algebraic Combinatorics Research Group, H--1117 Budapest, P\'azm\'any P.\ s\'et\'any 1/C, 
Hungary and Alfr\'ed R\'enyi Institute of Mathematics, Hungarian Academy of Sciences. 
Email: patkosb@cs.elte.hu and patkos@renyi.hu. Research supported by the 
J\'anos Bolyai Research Scholarship of the Hungarian Academy of Sciences.} 
\and Fiona Skerman \thanks{University of Oxford, Department of Statistics, 1 South
Parks Road, Oxford OX1 3TG, United Kingdom. Email:
skerman@stats.ox.ac.uk}}
\begin{document}
\maketitle

\begin{abstract}
In this paper, we study $(1 : b)$ Avoider-Enforcer games played on
the edge set of the complete graph on $n$ vertices. For every
constant $k\geq 3$ we analyse the $k$-star game, where Avoider tries
to avoid claiming $k$ edges incident to the same vertex. We consider
both versions of Avoider-Enforcer games -- the strict and the
monotone -- and for each provide explicit winning strategies for
both players. We determine the order of magnitude of the threshold
biases $f^{mon}_\cF$, $f^-_\cF$ and $f^+_\cF$, where $\cF$ is the
hypergraph of the game.

\end{abstract}

\section{Introduction}

Let $a$ and $b$ be two positive integers, let $X$ be a finite set
and let $\cF \subseteq 2^X$ be a family of subsets of $X$. In an
$(a:b)$ Avoider-Enforcer game $\cF$, two players, called Avoider and
Enforcer, alternately claim $a$ and $b$ previously unclaimed
elements of $X$ per move, respectively. If the number of unclaimed
elements is strictly less than $a$ (respectively $b$) before
Avoider's (respectively Enforcer's) move, then he claims all these
elements. The game ends when all the elements of $X$ have been
claimed by either of the players. Avoider loses the game if by the
end of the game he has claimed all the elements of some $F\in \cF$,
and wins otherwise. Throughout this paper we assume that Avoider is
the first player to play, although usually it makes very little
difference. We refer to $X$ as the \emph{board} of the game, to
$\cF$ as the \emph{target sets}, and to $a$ and $b$ as the
\emph{bias} of Avoider and Enforcer, respectively. Since the pair
$(X,\cF)$ is a hypergraph that represents the game, we often refer
to $\cF$ as the hypergraph of the game, or as the game itself.

Avoider-Enforcer games are the mis\`{e}re version of the well-studied Maker-Breaker games. In an $(a:b)$ Maker-Breaker game $\cF$,
the two players are called Maker and Breaker, they claim
respectively $a$ and $b$ elements of $X$ per move, and Maker wins if
and only if by the end of the game he has claimed all the elements
of some $F\in \cF$. Both Maker-Breaker and Avoider-Enforcer games
are finite, perfect information games, and there is no possibility
of a draw. Hence, for every given setup -- $a,b,\cF$ -- one of the
players has a winning strategy. We say that this player wins the
game.

It is very natural to play both Avoider-Enforcer and Maker-Breaker
games on the edge set of a given graph $G$, and specifically for
$G=K_n$, the complete graph on $n$ vertices. In this case the board
is $X=E(K_n)$ and the target sets are $\cF\subseteq 2^{E(K_n)}$. For
example: in the \emph{connectivity} game ${\cC_n}$ the target sets
are all edge sets of connected graphs on $n$ vertices;
in the \emph{perfect matching} game ${\cM_n}$ the target sets are
all graphs on $n$ vertices containing a perfect matching (we assume
$n$ is even here); in the \emph{Hamiltonicity} game ${\cH_n}$ the
target sets are all edge sets of graphs on $n$ vertices containing a
Hamilton cycle. We usually omit the subindex $n$ in our notation.
These three games were initially studied in Maker-Breaker version by
Chv\'{a}tal and Erd\H{o}s in their seminal paper \cite{CE}.

Many natural games played on the edges of $K_n$ (including all the
above mentioned ones) are drastically in favor of Maker, i.e.\ Maker
wins in the unbiased $(1:1)$ version in (almost) minimal number of moves required to create a winning set. Therefore, it makes
sense to give more power to Breaker in order to even out the odds,
and typically the $(1:b)$ version is considered. In addition,
Maker-Breaker games are \emph{bias monotone}: if Maker wins some
game $\cF$ with bias $(a:b)$, he also wins this game with bias
$(a':b')$, for every $a' \geq a$ and $b' \leq b$. This bias
monotonicity enables the definition of the \emph{threshold bias}:
for a given hypergraph $\cF$, the threshold bias $f_\cF$ is the
unique integer for which Maker wins the $(1:b)$ game $\cF$ for every
$b < f_\cF$, and Breaker wins the $(1:b)$ game $\cF$ for every $b
\ge f_\cF$.

Unfortunately, Avoider-Enforcer games are not bias monotone
in general (see e.g.\ \cite{HKSS10}, \cite{HKS07}): although
intuitively each player wishes to claim as few elements as possible,
it is sometimes a disadvantage to claim fewer elements per move, for
any of the players. This makes the analysis of these games much more
difficult, and it is not possible to define the threshold bias in
the same manner as in Maker-Breaker games. Therefore, Hefetz,
Krivelevich and Szab\'o introduced in \cite{HKS07} the following
parameters. The \emph{lower threshold bias} $f^-_\cF$ is the largest
integer such that Enforcer wins the $(1:b)$ game $\cF$ for every
$b\leq f^-_\cF$. The \emph{upper threshold bias} $f^+_\cF$ is the
smallest non-negative integer such that Avoider wins the $(1:b)$
game $\cF$ for every $b>f^+_\cF$. Except for some trivial cases,
$f^-_\cF$ and $f^+_\cF$ always exist and satisfy $f^-_\cF \leq
f^+_\cF$. When $f^-_\cF = f^+_\cF$ we call this number $f_\cF$ and
refer to it as the \textit{threshold bias} of the game $\cF$.

In order to overcome this bias monotonicity obstacle, Hefetz, Krivelevich,
Stojakovi\'c and Szab\'o proposed in \cite{HKSS10} a bias monotone
version for Avoider-Enforcer games: they suggested that Avoider and
Enforcer will claim \textbf{at least} $a$ and $b$ board elements per
move, respectively. It is easy to see that this new version is
indeed bias monotone, i.e.\ each player can only benefit from
lowering his bias. This fact allowed them to define for any given
hypergraph $\cF$ the \emph{monotone threshold bias} $f^{mon}_\cF$ as
the largest non-negative integer for which Enforcer wins the $(1:b)$
game $\cF$ under the new set of rules if and only if $b\leq
f^{mon}_\cF$. Throughout this paper we refer to this new set of
rules as the \emph{monotone} rules, to distinguish it from the
\emph{strict} rules. Accordingly, we refer to the games played under
the two sets of rules as monotone games and as strict games, respectively.

Interestingly, these seemingly minor adjustments in the rules can
completely change the outcome of the game. For example, even in such
a natural game as the connectivity game, the two versions of the
game are essentially different. In \cite{HKS07} it was shown that
Avoider wins the strict $(1:b)$ connectivity game played on $E(K_n)$
if and only if at the end of the game he has at most $n-2$ edges,
therefore the threshold bias exists and is of linear order. On the
other hand, the monotone threshold bias for this game is of order
$\frac {n}{\ln n}$ \cite{HKSS10,KS}.

Naturally, one may ask about the relationship between $f_{\cF}^{-}$,
$f_{\cF}^{+}$ and $f_{\cF}^{mon}$. Specifically, it could be
expected that $f_{\cF}^{-}\leq f_{\cF}^{mon} \leq f_{\cF}^+$ holds
for every family $\cF$. The above mentioned connectivity game shows
that this is not true in general, even when there exists a threshold
bias in the strict game.

In \cite{HKS07}, Hefetz, Krivelevich and Szab\'o provided a general
sufficient condition for Avoider's win in $(a:b)$ Avoider-Enforcer
games played under both sets of rules. This criterion takes only
Avoider's bias into account. In \cite{MBed}, Bednarska-Bzd\c{e}ga
introduced a new sufficient condition for Avoider's win under both
sets of rules, which depends on both parameters $a$ and $b$, and
gives a better result than the one in~\cite{HKS07} in cases where
the hypergraph of the game has rank smaller than $b$.

In \cite{HKSS10}, Hefetz et al.~investigated $(1:b)$
Avoider-Enforcer games played on the edge set of $K_n$, where
Avoider wants to avoid claiming a copy of some fixed graph $H$. In
this case $X=E(K_n)$, and $\cF=\cK_H \subseteq 2^{E(K_n)}$ consists
of all subgraphs of $K_n$ containing $H$ as a subgraph. These games
are referred to as $H$-games. They conjectured that for any fixed
graph $H$, the thresholds $f^-_{\cK_H}$ and $f^+_{\cK_H}$ are not of
the same order of magnitude, and wondered about the connection between monotone
$H$-games and strict $H^-$-games, where $H^-$ is $H$ with one edge
missing. They investigated $H$-games where $H=K_3$ (a triangle) and
$H=P_3$ (a path on three vertices) and established the following:
$$f^{mon}_{\cK_{P_3}}=\binom{n}{2}-\left\lfloor\frac
n2\right\rfloor-1,\quad f^+_{\cK_{P_3}}=\binom{n}{2}-2,\quad
f^-_{\cK_{P_3}}=\Theta(n^\frac 32)\enspace \mbox{ and } \enspace
f^{mon}_{\cK_{K_3}}=\Theta(n^\frac 32).$$ They used this example to
support their conjecture, as $f^+_{\cK_{P_3}}$ and $f^-_{\cK_{P_3}}$
are indeed not of the same order. They also noted that
$f^{mon}_{\cK_{K_3}}$ and $f^-_{\cK_{P_3}}$ are of the same order,
and that $P_3=K^-_3$. Bednarska-Bzd\c{e}ga established in
\cite{MBed} general upper and lower bounds on $f^+_{\cK_H}$,
$f^-_{\cK_H}$ and $f^{mon}_{\cK_H}$ for every fixed graph $H$, but these bounds are not tight for every graph $H$. In order to
prove our main result of this paper we prove a number theoretic
fact, which we later use independently in order to improve one of
Bednarska-Bzd\c{e}ga's bounds. We elaborate on that in \sref{apl}.

Our main objective in this paper is to study monotone and strict
$H$-games played on the edges of $K_n$, where $H$ is the $k$-star
$K_{1,k}$, denoted by $\cS_k$, for any fixed $k \geq 3$. We refer to
this game as the \emph{star game}, or more specifically, for a given
$k$, we call this game the \emph{$k$-star game}. Studying the star
game is very natural, since avoiding a $k$-star in Avoider's graph
is exactly keeping its maximal degree strictly below $k$. We
analyse this game, provide explicit winning strategies for both
players under both sets of rules, and obtain the following.

\begin{theorem}\label{thm:mainTheorem}
For every $k \ge 3$ and for every large enough $n$ the following
bounds hold:
\begin{enumerate}[(i)]
\item $\frac 25 n^{\frac{k}{k-1}} \le f^{mon}_{\cK_{\cS_k}} \le 2n^{\frac{k}{k-1}}$;
\item $\frac 15 n^{\frac{k}{k-1}} \le f^{ + }_{\cK_{\cS_k}} \le 2n^{\frac{k}{k-1}}$;
\item $\frac 12 n^{\frac{k+1}{k}} \le f^{ - }_{\cK_{\cS_k}} \le 2n^{\frac{k+1}{k}}$.
\end{enumerate}
\end{theorem}

These results show that $f^-_{\cK_{\cS_k}}$ and $f^+_{\cK_{\cS_k}}$
are not of the same order for any given $k \ge 3$, supporting the conjecture of Hefetz et al.~from~\cite{HKSS10}. In addition, as
$\cS^-_k = \cS_{k-1}$, an immediate consequence of
Theorem~\ref{thm:mainTheorem} is that $f^{mon}_{\cK_{\cS_k}}$ and
$f^-_{\cK_{\cS^-_k}}$ are of the same order, showing a strong
connection between the monotone $H$-game and the strict $H^-$-game
in this case. Note that $\cS_2 = P_3$, so the $k$-star game for
$k=2$ is already covered in~\cite{HKSS10}. In fact, the results
there match ours, if we generalize Theorem~\ref{thm:mainTheorem} to
include the case $k=2$. However, since these results are known, and
in order to avoid some technical difficulties in our proofs, we only
consider the case $k \ge 3$. 

The outcome of some $(1:b)$ positional games played on the edges of
$K_n$, where the target sets possess some graph property $\cP$, is
the same as in the corresponding games where the players play
randomly. This phenomenon was first observed by Chvat\'{a}l and
Erd\H{o}s in \cite{CE} for the Maker-Breaker connectivity game, and
is known as the \emph{random graph intuition}. The reason for this name is that
when both players play randomly the $(1:b)$ game, the graph of the
player with bias $1$ (either Maker or Avoider) at the end of the
game satisfies $G \sim G(n,m)$, where $m=\lceil \frac{1}{b+1} {n \choose 2}\rceil$. For this given $m$, the graph $G(n,m)$ behaves in many ways similarly to $G(n,\frac{1}{b+1})$, the random graph
on $n$ vertices where each potential edge appears in the graph
independently with probability $\frac{1}{b+1}$ \cite{JLR}. In other words, the
threshold bias $b^*$ for these games is asymptotically equal to $1 /
p^*$, where $p^*$ is the threshold probability for the appearance of
$\cP$ in $G \sim G(n,p)$.

The $k$-star game is a very good example for this phenomenon, as
indeed the properties of the random graph $G \sim G(n,\frac 1{b+1})$
suggest the outcome of the Avoider-Enforcer $(1:b)$ $k$-star game.
All the following statements about $G$ hold w.h.p.\ (i.e.\ with
probability tending to 1 as $n$ tends to infinity). For details the
reader may refer to \cite{Boll}, Theorem 3.1.
\begin{itemize}
\item For $b = \omega(n^{\frac{k}{k-1}})$ the maximal degree in $G$
is at most $k-2$, and Avoider wins the $(1:b)$ game (both strict
and monotone).

\item At $b = \Theta(n^{\frac{k}{k-1}})$, vertices of degree $k-1$
emerge in $G$. If Avoider claims the last edge in the $(1:b)$ game,
the appearance of a vertex of degree $k-1$ in his graph before the
last round means he loses, and this is indeed the order of magnitude
of $f^{mon}_{\cS_k}$ and $f^+_{\cS_k}$, where presumably Avoider
claims the last edge.

\item When $b=\omega\left(n^{\frac{k+1}{k}}\right)$ and $ b=o\left(n^{\frac{k}{k-1}}\right)$, the maximal degree in $G$ is exactly $k-1$. The
outcome of the strict $(1:b)$ game heavily depends on the number of
free edges Avoider will be able to choose from in his last move, and
so the outcome oscillates.

\item Finally, for $b \le Cn^{\frac{k+1}{k}}$, where $C$ is a
sufficiently small constant, vertices of degree $k$ emerge in $G$,
and Enforcer wins the $(1:b)$ game (both strict and monotone).
\end{itemize}

The rest of the paper is organized as follows: in \sref{avoider} we
provide Avoider's strategy for the $k$-star game which applies for
both versions of the game. In \sref{enfStar} we provide Enforcer's
strategies for the $k$-star game, one strategy for the monotone game
and one for the strict game. In~\sref{apl} we improve one of Bednarska-Bzd\c{e}ga's bounds for general $H$-games. Finally, in~\sref{rmrks} we present
some concluding remarks and open problems.

\section{Preliminaries}
\label{sec:prel}

Throughout this paper we use the following notation.

A previously unclaimed edge is called a \emph{free} edge. The act of
claiming one free edge by one of the players is called a
\emph{step}. In the strict game, Enforcer's $b$ (respectively
Avoider's $1$) successive steps are called a \emph{move}. In the
monotone game, each move consists of at least $b$ steps,
respectively at least one step. A \emph{round} in the game consists
of one move of the first player (Avoider), followed by one move of
the second player (Enforcer). Whenever one of the players claims an
edge incident to some vertex $u$, we say that the player
\emph{touched} $u$.

Our graph-theoretic notation is standard and follows that of
\cite{West}. In particular, throughout the paper $G$ stands for a
simple graph with vertex set $V=V(G)$ and edge set $E=E(G)$. For any
subset $U \subseteq V$ we say that an edge $uv$ lies \emph{inside}
$U$ if $u,v \in U$. For $i\geq 0$, we denote by $A_i$ and $E_i$ the graphs with
vertex set $V$, whose edges were claimed by Avoider, respectively
Enforcer, in the first $i$ rounds. For every vertex $v\in V$ and
every $i \ge 0$, let $d_{A_i}(v)$ and $d_{E_i}(v)$ denote the degree
of $v$ in $A_i$, respectively $E_i$. We sometimes omit the subindex
$i$ when its value is clear or irrelevant. In these cases we also
refer to $d_{A}(v)$ as the $A$-degree of $v$. Whenever we consider
the end of the $i$th round for the case $i=0$, we simply refer to
the beginning of the game, before any move was played.

The set of all free edges at the end of the $i$th round is denoted
by $F_i$. A free edge is called a \emph{threat} if it is incident to
a vertex of $A$-degree $k-1$.

For the sake of simplicity and clarity of presentation, no real
effort has been made here to optimize the constants appearing in our
results. We also omit floor and ceiling signs whenever these are not
crucial. Our results are asymptotic in nature and whenever necessary
we assume that $n$ is sufficiently large. We use $o(1)$ to denote a 
positive function of $n$, tending to zero as $n$ tends to infinity.

For every two integers $n$ and $b$ let $r=r(n,b)$ be the integer for
which $1 \le r \le b+1$ and $\binom{n}{2}\equiv r$ mod $(b+1)$ hold.
The value of $r$ is the number of free edges before the last round
of the strict game, and since Avoider is the first player, $r$ is
actually the number of edges which remain for Avoider to choose from
in his last move. Therefore, this value may be very significant in
determining the identity of the winner in the strict game. In order
to estimate $r$ in some cases, we need the following two number
theoretical statements.

\begin{fact}\label{fact:largeRemainder}
Let $c$ and $\alpha$ be two constants such that either $\alpha = 1$
and $c \ge 1$, or $\alpha \in (1,2)$ and $c > 0$. For any
sufficiently large integer $n$ there exists an integer $q =
\left(2-o(1)\right)cn^\alpha$ such that the remainder of the
division of $\binom n2$ by $q$ is larger than $cn^\alpha$.
\end{fact}

\begin{proof}
Let $N = \binom n2$, $N' = N - cn^\alpha - 1$, $m = \lceil\frac
1{4c} n^{2-\alpha}\rceil$ and $q = \lfloor \frac {N'}m \rfloor$.
Note that $N' = qm + r$ for some $0 \le r < m$ and that $q =
\left(2-o(1)\right)cn^\alpha$. Since $N-qm = cn^\alpha + r + 1 < q$,
it follows that the remainder of the division of $N$ by $q$ is
larger than $cn^\alpha$.
\end{proof}

\begin{fact}\label{fact:smallRemainder}
For every sufficiently large integer $n$ and for every constant $k
\ge 3$ there exists a $c = c(n,k)$ such that $\frac 15 < c < \frac
14$ and the remainder $r$ of the division of $\binom n2$ by
$cn^{\frac{k}{k-1}}$ is positive and satisfies $r = o(n)$.
\end{fact}

\begin{proof}
Let $N={n \choose 2}$ and let $M = \{2\lceil
n^{1-\frac{1}{k-1}}\rceil+i: 0\leq i \leq 4\}$. The least common
multiple of the elements in $M$ is 
 at least $n^{2.5}$ so
there exists an element $m \in M$ which does not divide $N$. Let $q
= \lfloor\frac Nm\rfloor$ and note that $q = \left(\frac 14
-o(1)\right) n^{\frac{k}{k-1}}$. Since $N-m < qm < N$, the remainder
$r$ of the division of $N$ by $q$ satisfies $0 < r < m = o(n)$.
\end{proof}


\section{Avoider's strategy}
\label{sec:avoider}

In this section we establish upper bounds on the threshold biases
$f^{mon}_{\cK_{\cS_k}}$, $f^+_{\cK_{\cS_k}}$ and
$f^-_{\cK_{\cS_k}}$.

We provide Avoider with the following trivial strategy $\cS_A$: in
every move Avoider claims one arbitrary edge which does not increase
the maximal degree in his graph if such an edge exists, and an
arbitrary edge otherwise. Clearly Avoider can follow this strategy.
Note that this is a valid strategy for both the monotone and the
strict versions of the game.

Consider the course of a game (either strict or monotone) in which
Avoider plays according to $\cS_A$ and Enforcer plays according to
some fixed strategy. For every $i$, let $I_i$ denote the set of
vertices of maximal $A$-degree at the end of round $i$. Let $s$ be
the maximal $A$-degree at the end of the game. For every $0 \le j
\le s$, let $i_j$ be the largest integer such that the maximal
$A$-degree at the end of round $i_j$ is $j$. Note that the maximal
$A$-degree is never increased by more than one according to $\cS_A$,
and so $0 = i_0 < i_1 < \cdots < i_s$.

\begin{lemma}\label{lem:stages}
$|F_{i_j}| \le (2^{j-1}+o(1))\frac{n^{j+2}}{b^j}$ for every $0 \le j
\le s$.
\end{lemma}

\begin{proof}
Observe that at the end of round $i_j$ (for every $j$) every free
edge has at least one endpoint in $I_{i_j}$, as otherwise Avoider
will not increase the maximal degree in his graph in his subsequent
move. Therefore, if $|I_{i_j}| \le
(2^{j-1}+o(1))\frac{n^{j+1}}{b^j}$ then $|F_{i_j}|\le |I_{i_j}|\cdot n \le
(2^{j-1}+o(1))\frac{n^{j+2}}{b^j}$. Since the number of free edges
at the beginning of the game is obviously $\binom n2 \le 2^{-1}n^2$,
it suffices to show that if $|F_{i_j}| \le
(2^{j-1}+o(1))\frac{n^{j+2}}{b^j}$, then $|I_{i_{j+1}}| \le
(2^{j}+o(1))\frac{n^{j+2}}{b^{j+1}}$, for every $0 \le j < s$.

Indeed, as both players claim altogether at least $b+1$ edges, for every $0 \le j < s$ the number of rounds in the game
after round $i_j$ cannot be greater than
$\lceil|F_{i_j}|/(b+1)\rceil \le
(2^{j-1}+o(1))\frac{n^{j+2}}{b^{j+1}}$. Since Avoider claims exactly
one edge per move, in each round after round $i_j$ at most two new
vertices of $A$-degree $j+1$ appear. Therefore, $|I_{i_{j+1}}| \le
(2^{j}+o(1))\frac{n^{j+2}}{b^{j+1}}$.
\end{proof}

Now it is easy to see that $\cS_A$ is a winning strategy for Avoider
in the $(1:b)$ Avoider-Enforcer $k$-star game for any $b\geq
2n^{\frac{k}{k-1}}$ and under both sets of rules, thus obtaining the
upper bounds in Theorem~\ref{thm:mainTheorem}~$(i)$ and~$(ii)$.
Indeed, it follows by Lemma~\ref{lem:stages} that if $i_{k-2}$
exists then $|F_{i_{k-2}}| \le n^{\frac{k}{k-1}} < b$. Therefore no
threat appears before Avoider's last move and so he wins.

We now prove the upper bound for $f^-_{\cK_{\cS_k}}$ given in
Theorem~\ref{thm:mainTheorem}~$(iii)$. By
Fact~\ref{fact:largeRemainder}, with $c = 1$ and $\alpha = \frac
{k+1}{k}$, there exists an integer $b = \left(2-o(1)\right)n^\frac
{k+1}{k}$ such that $r(n,b) > n^\frac {k+1}{k}$. Assume that $s \ge
k-1$ for this $b$ (otherwise Avoider obviously wins). It follows by
Lemma~\ref{lem:stages} that $|F_{i_{k-1}}| < n^{\frac{k+1}{k}}$.
However, by the assumption on $r$ there are more free edges than
$n^{\frac{k+1}{k}}$ before any move of Avoider, so $i_{k-1}$ is the last round of
the game, meaning $s = k-1$.

\begin{remark}\label{rem:kEqualsTwo}
All arguments in this section are still valid even if we include the
case $k=2$. 
\end{remark}

\section{Enforcer's strategies}
\label{sec:enfStar}
In this section we establish the lower bounds given in
Theorem~\ref{thm:mainTheorem}. Unlike Avoider's strategy, which was
valid for both versions of the game, here we distinguish between the
two cases. We start with the monotone game which is simpler to
analyse and establish the lower bound on $f^{mon}_{\cK_{\cS_k}}$.
Then we proceed to the strict game, explain the adjustments we make
to Enforcer's strategy and establish the lower bounds on
$f^+_{\cK_{\cS_k}}$ and $f^-_{\cK_{\cS_k}}$.

\subsection{The monotone game}\label{sec:enfStarMon}
We provide a strategy for Enforcer for the monotone $(1:b)$ $k$-star
game for $b = \frac25 n^{\frac{k}{k-1}}$. At any point during the
game, let $I$ denote the set of isolated vertices in Enforcer's
graph, and let $C = V \setminus I$. Furthermore, let $I_i$ and $C_i$
denote the respective sets of vertices at the end of the $i$th
round. Initially, of course, $I_0=V$ and $C_0=\emptyset$. Whenever
Enforcer touches a vertex previously isolated in his graph, we say
that he moved that vertex from $I$ to $C$.

For every $i \geq 0$, Enforcer plays his $(i+1)$st move as follows.
\begin{description}
\item [(1)]If there exists a vertex of $A$-degree at least $k$, or
if there are at most $b$ free edges remaining, Enforcer claims all
free edges on the board. We refer to this move as the \emph{trivial
move}.
\item [(2)] Otherwise, if there exists a vertex $v \in I$ of $A$-degree $k-1$,
then Enforcer claims all free edges on the board, except one,
incident to $v$. We refer to this move as the \emph{end move}.
\item [(3)] Otherwise, let $v_1^{(i)},\ldots,
v_{|I_i|}^{(i)}$ be an enumeration of the vertices in $I_i$ such
that for every $1\leq j <|I_i|$, $d_{A_{i+1}}(v_j^{(i)})\leq
d_{A_{i+1}}(v_{j+1}^{(i)})$. Let $I_{i,j} = \{v_1^{(i)},\ldots,
v_j^{(i)}\}$ and let $s_i$ be the smallest integer such that the
number of free edges inside $C_i \cup I_{i,s_i}$ is at least $b$.
Enforcer claims all the free edges inside $C_i \cup I_{i,s_i}$. We refer to this move as the \emph{base move}.
\end{description}

We have to show that Enforcer can follow the proposed strategy, and
that by doing so he wins the game. Starting with the former, it is
evident that Enforcer can play the trivial move; Enforcer can play
the end move since there are more than $b$ free edges on the board,
and since there are $n-k$ free edges incident to $v$; finally,
Enforcer can play the base move since $C_i \cup I_{i,s_i} = V$ for
$s_i = |I_i|$ and there are more than $b$ free edges on the board.

We now prove that the proposed strategy is indeed a winning strategy
for Enforcer. Consider the course of the game in which Enforcer
plays according to the proposed strategy and Avoider plays according
to some fixed arbitrary strategy. If at any point during the game
the maximal $A$-degree in the graph increases to at least $k$ then
Enforcer wins. He also wins if he plays the end move at some point.
So assume for contradiction that neither of these events happen.
Therefore, by the description of his strategy, it is clear that
Enforcer plays the base move for $l$ rounds, for some $l \ge 0$, and
then either the game ends or in his last move he plays the trivial
move since there are at  most $b$ free edges remaining.

\begin{observation}\label{obs:monotoneEnforcer}
Throughout the game, the following properties hold.
\begin{enumerate}[(i)]
\item There are at least $n-k$ free edges incident to every vertex in $I$.
\item After every move, by either player, every
free edge has at least one endpoint in $I$.
\item The number of edges claimed by both players in each round of
the game is at most $(1+o(1))b$.
\end{enumerate}
\end{observation}

\emph{Proof.}
\begin{enumerate}[$(i)$]
\item This is obvious since every vertex in $I$ is isolated in Enforcer's
graph and every vertex has $A$-degree less than $k$.
\item The claim is true after each base move played by Enforcer by his
strategy, and there are no free edges left after he plays the
trivial move, if he does. Recall that by assumption he never plays
the end move. In addition, Avoider does not change the set $I$, so
the claim remains true after his moves as well.
\item Whenever Enforcer plays the base move he does not
claim more than $b+n=(1+o(1))b$ edges. If he plays the trivial move
this is obviously still true. Finally, Avoider claims at most
$\frac{kn}{2} = o(b)$ edges throughout the game, otherwise a vertex
of $A$-degree at least $k$ must exist.{\hfill $\Box$\medskip}
\end{enumerate}

Now we wish to estimate the $A$-degrees of vertices in $I$. Let
$T(i):=\frac{1}{|I_i|}\sum_{v\in I_i}d_{A_i}(v)$ denote the average
$A$-degree of the vertices in $I_i$ at the end of round $i$. Note
that by definition $T(0)=0$, and that $T(i) > T(i-1)$ if Enforcer
plays the base move in his $i$th move. Indeed, since throughout the
game all free edges have at least one endpoint in $I$, Avoider in
his $i$th move increases the sum of $A$-degrees in $I$ (while not
changing the set itself), and Enforcer in his subsequent move
removes from $I$ vertices of minimal $A$-degree, so he does not decrease
the average $A$-degree in $I$.

\begin{claim}\label{clm:avgDegreeMonotone}
For every $0\leq j \leq k-2$, the following holds. If $0 < |I_i| <
\frac{9}{10}n^{1-\frac{j}{k-1}}$ for some $i$, then $T(i)\geq j$.
\end{claim}

\begin{proof}
We prove the claim by induction on $j$. The claim trivially holds
for $j=0$, as $T(i) \ge 0$ for every $i$. Suppose now for
contradiction that for some $1 \le j \le k-2$, the claim holds for
$j-1$, but not for $j$. Then there exists an integer $i$ such that
$0 < |I_i| < \frac{9}{10} n^{1-\frac{j}{k-1}}$, but $T(i) < j$,
which implies $\sum_{v \in I_i}d_{A_i}(v) <j|I_i|$. Let $i_0 \le i$
be the minimal index such that $T(i_0)\geq j-1$ (by the induction
hypothesis and the size of $I_i$, such an index exists), and for
every $i_0 \leq s \leq i$ let $W(s):=\sum_{v \in
I_s}\left(d_{A_s}(v)-(j-1)\right)$. Note that $W(i) < |I_i|$ by the
assumption on the index $i$.

Since $T(i)<j$, there are vertices of $A$-degree less than $j$ in
$I_i$, and therefore, according to his strategy, Enforcer has only
moved vertices of $A$-degree less than $j$ from $I$ to $C$ in his
first $i$ moves. In addition, Avoider increases the sum of
$A$-degrees of the vertices in $I$ in each of his moves. It follows
that $W(s+1) > W(s)$ for every $i_0 \le s < i$. Since $W(i_0) \ge 0$
by definition of $i_0$, and since $W(s)$ is an integer for every
$s$, we get that $i - i_0 + 1 \leq W(i) + 1 \le |I_i| < \frac{9}{10}
n^{1-\frac{j}{k-1}}$. It follows that between rounds $i_0$ and $i$,
including round $i_0$ if $i_0 > 0$, Enforcer has claimed at most
$(1+o(1))b|I_i| < \frac25 n^{2-\frac{j-1}{k-1}}$ edges.

On the other hand, consider the vertices that were moved from $I$ to
$C$ by Enforcer between rounds $i_0$ and $i$. Let $I^* = I_{i_0-1}$
if $i_0 > 0$, and $I^* = I_0 = V$ otherwise (note that $i_0 = 0$ if
and only if $j = 1$). Since throughout the game every vertex in $I$
has at least $n-k$ free edges incident to it, and by using the
induction hypothesis, we conclude that during the specified rounds
Enforcer must have claimed at least
\begin{align*}
\frac{\left(|I^*| - |I_i|\right) \left(|C_i| - k\right)}{2} & \geq
\left(1-o\left(1\right)\right)\left(\frac{9}{10}
n^{1-\frac{j-1}{k-1}} - \frac{9}{10}
n^{1-\frac{j}{k-1}}\right)\frac{|C_i|}{2}\\
& = \left(1-o\left(1\right)\right)\frac{9}{20}
n^{2-\frac{j-1}{k-1}} \\
&> \frac25 n^{2-\frac{j-1}{k-1}}
\end{align*}
edges, a
contradiction.
\end{proof}

Let $i$ be the maximal index such that $|I_i| > 0$ and $T(i)< k-2$
(there exists such an index since both inequalities hold for $i=0$).
By Claim~\ref{clm:avgDegreeMonotone} we get $|I_{i}|\geq
\frac{9}{10} n^{\frac{1}{k-1}}$. Hence, either $|I_i|\geq \frac{n}{1000}$ 
and then $|F_i| = \Theta(n^2)=\omega(b)$, or
$$|F_i| \ge |I_i|(n-k) - \binom{|I_i|}2 \ge \left(1-\frac{1}{100}\right)|I_i|n >
\frac{89}{100}n^{1+\frac{1}{k-1}}> \frac{22}{10}b.$$

Therefore, $|F_{i+1}| > \frac{11}{10}b$ by Part $(iii)$ of
Observation~\ref{obs:monotoneEnforcer}, and $T(i+1) \ge k-2$ by
definition of $i$. Hence, after Avoider's $(i+2)$nd move, either
there exists a vertex of $A$-degree at least $k$, or there exists a
vertex $v \in I$ with $A$-degree $k-1$, while there are still more
than $b$ free edges on the board, in which case Enforcer plays the
end move. In either case, this is a contradiction to the assumption
on Avoider's strategy. This completes the proof.

\subsection{The strict game}\label{sec:enfStarStrict}
Recall that $r = r(n,b)$ denotes the integer which satisfies $1\le
r\le b+1$ and $\binom{n}{2}\equiv r$ mod $(b+1)$, i.e.\ the number of
free edges at the beginning of the last round of the game. Let

\begin{itemize}
\item[] $b^+_{n,k}=\max\left\{b \le \frac 14 n^{\frac{k}{k-1}}: r(n,b) \le
\frac 58 \frac{n^{k+1}}{(2b)^{k-1}}\right\}, \mbox{ and}$
\item[] $b^-_{n,k}=\max\left\{b \leq \frac 14 n^{\frac{k}{k-1}}: r(n,b') \le
\frac 58 \frac{n^{k+1}}{(2b')^{k-1}}\textrm{ for every } 1 \le b'
\le b\right\}.$
\end{itemize}

\begin{claim}\label{clm:lowerBounds}
For every sufficiently large integer $n$ and for every integer $k
\ge 3$ the following bounds hold:
\begin{enumerate}[$(i)$]
\item $b^+_{n,k} \ge \frac{1}{5}n^{\frac{k}{k-1}}$;
\item $b^-_{n,k} \ge \frac{1}{2}n^{\frac{k+1}{k}}$.
\end{enumerate}
\end{claim}

\emph{Proof.}
\begin{enumerate}[$(i)$]
\item By Fact~\ref{fact:smallRemainder} there exists an integer
$\frac{1}{5}n^{\frac{k}{k-1}} \le b \le
\frac{1}{4}n^{\frac{k}{k-1}}$ such that $r(n,b) = o(n)$, and since
$\frac{5}{8}\frac{n^{k+1}}{(2b)^{k-1}}=\Theta(n)$ in this case, the
desired inequality holds.
\item Note that $$b \le \frac{1}{2} n^{\frac {k+1}{k}} \Longrightarrow
(2b)^k \le n^{k+1} \Longrightarrow b \le \frac{1}{2} \frac
{n^{k+1}}{(2b)^{k-1}},$$ and since $r \le b+1$ trivially holds, we
get $b^-_{n,k} \ge \frac{1}{2} n^{\frac {k+1}{k}}$.
{\hfill$\Box$\medskip}
\end{enumerate}

The lower bounds in Theorem~\ref{thm:mainTheorem} $(ii)$ and~$(iii)$
follow directly from Claim~\ref{clm:lowerBounds} and the following lemma.

\begin{lemma}\label{lem:bounds}
$b^+_{n,k}\leq f^+_{\cK_{\cS_k}}$ and $b^-_{n,k}\leq
f^-_{\cK_{\cS_k}}$ hold for every $k \ge 3$ and sufficiently large
$n$.
\end{lemma}

\begin{proof}
Throughout this proof we assume that Enforcer's bias $b$ satisfies
$b\leq \frac 14 n^{\frac{k}{k-1}}$. For simplicity, we first assume
that $b$ also satisfies $b=\omega(n)$. We propose a strategy for
Enforcer which is very similar to the proposed strategy in the
monotone game. However, some modifications are inevitable. One major
difference between the two versions of the game is that the
appearance of one threat (recall that a threat is a free edge
incident to a vertex of $A$-degree $k-1$) does not secure Enforcer's
win, so he has to make sure that $r$ threats appear before the last
round. We therefore say that the game is in a \emph{winning
position} if either the maximal degree in Avoider's graph is at
least $k$ or there exist at least $r$ threats. Since Enforcer cannot
increase Avoider's degrees or the number of threats, Enforcer wins
the game if and only if the game is in a winning position after
Avoider's penultimate move. For convenience we denote $\ell =
\left\lceil\binom n2/(b+1)\right\rceil - 1$ (i.e.\ the game lasts
exactly $\ell+1$ rounds).

Another difference between the two versions of the game is that in
the strict game Enforcer cannot maintain the property that every
free edge is incident to at least one vertex isolated in his graph.
However, he is able to maintain a partition $V=I \cup C$ (where
$I_i$ and $C_i$ denote the respective sets at the end of the $i$th
round) with some similar properties. The exact construction of the
sets $I$ and $C$ will be explained shortly. Initially, as in the
monotone game, $I_0=V$ and $C_0=\emptyset$. Once again we denote by
$T(i):=\frac{1}{|I_i|}\sum_{v\in I_i}d_{A_i}(v)$ the average
$A$-degree of the vertices in $I_i$ at the end of round $i$.

Enforcer's strategy involves dividing the course of the game into
two stages. The game begins at Stage I; for every $0 < i < \ell$, if
the game is in a winning position before Enforcer's $i$th move then
Stage~I is over and Enforcer immediately proceeds to Stage~II.
Otherwise, he keeps playing in Stage~I. If before Enforcer's
$\ell$th move the game is still in Stage~I, he proceeds to Stage~II
even if the game is not in a winning position. So, for some $1 \le i
\le \ell$ Avoider's $i$th move is the last move in Stage~I and
Enforcer's $i$th move is the first move in Stage~II. In each stage,
Enforcer plays as follows.

\begin{description}
\item [Stage~I:] For every $i \geq 0$ such that Enforcer plays his
$(i+1)$st move in this stage, let $v_1^{(i)},\ldots,
v_{|I_i|}^{(i)}$ be an enumeration of the vertices in $I_i$ for
which $d_{A_{i+1}}(v_j^{(i)})\leq d_{A_{i+1}}(v_{j+1}^{(i)})$, and
if $d_{A_{i+1}}(v_j^{(i)}) = d_{A_{i+1}}(v_{j+1}^{(i)})$ then
$d_{E_{i}}(v_j^{(i)})\geq d_{E_{i}}(v_{j+1}^{(i)})$, for $1\leq j
<|I_i|$. Let $I_{i,j} = \{v_1^{(i)},\ldots, v_j^{(i)}\}$ and let
$s_i$ be the largest integer such that the number of free edges
inside $C_i \cup I_{i,s_i}$ is at most $b$. Every move consists of
two parts.

In the first part of every move Enforcer claims all the free edges
inside $C_i \cup I_{i,s_i}$ and he moves $I_{i,s_i}$ from $I$ to $C$,
i.e.\ defines $C_{i+1}:=C_i \cup I_{i,s_i}$ and $I_{i+1}:=I_i
\setminus I_{i,s_i}$.

For the second part of every move, let $l_{i+1}$ denote the number
of edges Enforcer must claim in order to complete his $(i+1)$st
move. For every vertex $v \in
\{v^{(i)}_{s_i+1},\dots,v^{(i)}_{s_i+4k}$\}, Enforcer claims either
$\lfloor \frac{l_{i+1}}{4k} \rfloor$ or $\lceil \frac{l_{i+1}}{4k}
\rceil$ arbitrary free edges $vu$ such that $u \in C_{i+1}$, to get
a total of $l_{i+1}$ edges, thus completing his move. We say that
these edges are \emph{attached} to $v$.
\item [Stage~II:] In every step of every move in this stage, Enforcer
claims an arbitrary edge which is not a threat if such an edge
exists, and an arbitrary threat otherwise. He no longer maintains
the partition $V=I \cup C$.
\end{description}

First we show that Enforcer can follow the proposed strategy. This
is obvious for Stage~II and for the first part of every move in
Stage~I. Assume now that Enforcer is trying to play the second part
of his $i$th move in Stage~I for some $i > 0$, after playing
successfully all his previous moves according to the proposed
strategy, including the first part of the $i$th move. In particular,
the partition $V = I_i \cup C_i$ has been determined. It is easy to
see that at this point, exactly as in the monotone game, every free
edge has at least one endpoint in $I_{i}$. Therefore, if there are
less than $4k$ vertices in $I_i$ then there are only $O(n) = o(b)$
free edges remaining (by our assumption $b = \omega(n))$, which
implies $i \ge \ell$, in contradiction to the assumption that
Enforcer is playing his $i$th move in Stage~I.

Hence, it only remains to show that Enforcer will be able to attach
enough edges to every vertex among the first $4k$ of $I_i$. Observe
that $l_{i} < |C_{i}|$ by definition of $s_{i-1}$ and that $|C_j| =
\omega(1)$ for every $j > 0$. The following claim shows that
Enforcer can indeed follow the second part of his moves in Stage~I.

\begin{claim}\label{clm:extraEdges}
Throughout Stage~I there are at least $\left(\frac 34 - o(1)\right)|C|$ free
edges between every vertex in $I$ and $C$.
\end{claim}

\begin{proof}
Since $d_A(v) < k$ for every $v \in V$ throughout Stage~I, it
suffices to show that for every round $i$ in this stage, $d_{E_i}(v)
\le \left(\frac{1}{4} + o(1)\right)|C_{i}|$ for every $v \in I_i$.

Let $v \in I_i$ be a vertex that was touched by Enforcer in his
$i$th move. If Avoider does not touch $v$ in his $(i+1)$st move,
then in every proper enumeration of the vertices in $I_i$ before
Enforcer's $(i+1)$st move, $v$ will be among the first $4k$
vertices. Indeed, let $u$ be a vertex that was placed after the
first $s_i + 4k$ vertices of $I_{i-1}$ in the $i$th enumeration. By
the properties of the enumeration and our assumption we get
$d_{A_{i+1}}(v) = d_{A_{i}}(v) \le d_{A_{i}}(u) \le d_{A_{i+1}}(u)$.
In case of equality we get $d_{E_{i}}(v) > d_{E_{i-1}}(v) \ge
d_{E_{i-1}}(u) = d_{E_{i}}(u)$. Enforcer will then add $v$ to
$C_{i+1}$ since $b > 4kn$. So, $v$ remains in $I$ only if Avoider
touches it in his $(i+1)$st move and therefore every vertex can have
edges attached to it by Enforcer in at most $k$ rounds.

Now consider a vertex $v \in I_i$ for some $i>0$ (the claim is
trivial for $i=0$). Since $l_{j+1} < |C_{j+1}| \le |C_{i}|$ for
every $j < i$, the number of edges attached to $v$ cannot be more
than $k \lceil \frac{|C_{i}|}{4k}\rceil = \left(\frac{1}{4} +
o(1)\right)|C_{i}|$.
\end{proof}
We now wish to examine the course of the game in which Avoider plays
according to some fixed strategy and Enforcer plays according to the
proposed strategy, in order to obtain a sufficient condition for
Enforcer's win, thus proving the lemma. Note that if Enforcer plays
according to Stage~II of the strategy at any move before his
$\ell$th, he wins the game. Assume, then, that this does not happen.
It is immediate to observe that some properties hold exactly as in
the monotone game.
\begin{observation}\label{obs:strict}
The following properties hold throughout Stage~I.
\begin{enumerate}[(i)]
\item After every move, by either player, every free edge has at least one endpoint in $I$.
\item Enforcer has no edges inside $I$.
\item $T(0)=0$ and $T(i+1) > T(i)$.
\end{enumerate}
\end{observation}

The following claim is the strict analogue of
Claim~\ref{clm:avgDegreeMonotone}, showing that as $I$ gets smaller,
the average $A$-degree of its vertices becomes larger.

\begin{claim}\label{clm:avgDegreeStrict}
For every $0\leq j \leq k-2$ the following holds. If for some $i$
Enforcer plays his $i$th move according to Stage~I of his strategy
and $0 < |I_i| < \frac{9}{10} n\left(\frac {n}{2b}\right)^j$, then
$T(i)\geq j$.
\end{claim}

\begin{proof}
We prove the claim by induction on $j$. The claim trivially holds
for $j=0$, as $T(i) \ge 0$ for every $i$. Suppose now for
contradiction that for some $1 \le j \le k-2$, the claim holds for
$j-1$, but not for $j$. Then there exists an integer $i$ such that
$0 < |I_i| < \frac{9}{10} n\left(\frac {n}{2b}\right)^j$, but $T(i)
< j$.

As in the proof of Claim~\ref{clm:avgDegreeMonotone}, we denote by
$i_0 \le i$ the minimal index such that $T(i_0)\geq j-1$. Since the integer-valued weight function $W(s):=\sum_{v \in I_s}
\left(d_{A_s}(v)-(j-1)\right)$ is non-negative for $i_0$, and is
strictly increasing for $s \le i $, we conclude that $i - i_0 + 1\le
W(i) + 1 \le |I_i|$ and thus between rounds $i_0$ and $i$ Enforcer
has claimed at most $b|I_i| < \frac{9}{20} n^2 \left(\frac
{n}{2b}\right)^{j-1}$ edges.

We now show that according to his strategy Enforcer had to claim
more edges than that during these rounds. We distinguish between the
following cases:

For $j=1$, Avoider could not have claimed any edge inside $C_i$, so
the number of edges Enforcer had to claim is at least $\binom
{|C_i|}2 = (1-o(1))\binom n2 > \frac 9{20} n^2$.

For $j > 1$, note that $i_0 > 0$. If $|I_{i_0-1}| = \Theta(n)$ it
means that Enforcer claimed a quadratic number of edges in the
specified rounds (all edges inside $I_{i_0-1} \setminus I_i$, except
at most $kn$ edges that could have been claimed by Avoider). On the
other hand, if $|I_{i_0-1}| = o(n)$, note that Enforcer had to claim
all the edges between $I_{i_0-1} \setminus I_i$ and $C_{i_0-1}$ that
were free before round $i_0$, except at most $kn$ edges. By
Claim~\ref{clm:extraEdges} and the induction hypothesis, the number
of these edges is at least
$$\left(\frac{9}{10} n\left(\frac {n}{2b}\right)^{j-1} -
\frac{9}{10} n\left(\frac {n}{2b}\right)^j\right) \left(\frac 34
-o(1)\right) n
> \frac{9}{20} n^2 \left(\frac {n}{2b}\right)^{j-1}.$$
In either case, we get a contradiction.
\end{proof}

Let $g$ be the maximal index such that $|I_g| > 0$ and $T(g)< k-2$
(there exists such an index since both inequalities hold for $g=0$).
The next claim shows that as the game goes on after the $g$th round,
more and more vertices of $A$-degree $k-1$ appear in the graph.

\begin{claim}
\label{clm:finalstage} For every $i \ge 0$, after Avoider's
$(g+1+i)$th move either Avoider's graph contains an $\cS_k$ or there
are at least $i$ vertices in $I_{g+i}$ of $A$-degree $k-1$.
\end{claim}

\begin{proof}
Let $v_1,\dots,v_t$ denote the vertices of $I_{g+1}$ with $A$-degree
at most $k-3$ after the $(g+1)$st round, and let
$m=\sum_{i=1}^t\left(k-2-d_{A_{g+1}}(v_i)\right)$. If Avoider has
not yet created an $\cS_k$, all vertices have $A$-degree at most
$k-1$, thus after the $(g+1)$st round there are at least $m$
vertices in $I_{g+1}$ with $A$-degree $k-1$, as the average
$A$-degree in $I_{g+1}$ is at least $k-2$ by definition of $g$.
Since all free edges have at least one endpoint in $I$, as long as
Enforcer removes from $I$ only vertices of $A$-degree at most $k-2$,
in every round after the $(g+1+m)$th the number of vertices of
$A$-degree $k-1$ is increased, or an $\cS_k$ appears in Avoider's
graph. If Enforcer removes from $I$ a vertex of $A$-degree at least
$k-1$, then if the maximal degree in Avoider's graph is $k-1$ at
that point, it will be increased to $k$ in Avoider's subsequent move
(if such a move exists).
\end{proof}

By Claim~\ref{clm:avgDegreeStrict} we get $|I_{g}|\geq \frac{9}{10}
n\left(\frac {n}{2b}\right)^{k-2}$. Hence, either $I_g$ is of linear
order and then $|F_g| = \Theta(n^2)$ (since at most $k|I|$ edges can
be claimed inside $I$ throughout Stage~I), or $|C_g| = (1-o(1))n$
and then by Claim~\ref{clm:extraEdges} we get:
$$|F_g| \ge \left(\frac 34 - o(1)\right) |I_g|n \ge
\left(1-o(1)\right)\frac{27}{40} n^2\left(\frac
{n}{2b}\right)^{k-2}.$$ Thus, the number of rounds remaining in the
game after the $g$th round satisfies
\begin{equation}\label{eq:1}
\ell + 1 - g = \left\lceil\frac {|F_g|}{b+1}\right\rceil \ge
\left(1-o(1)\right)\frac{27}{20} n\left(\frac {n}{2b}\right)^{k-1}.
\end{equation}

Using our assumption $b\leq \frac 14 n^{\frac{k}{k-1}}$, a simple
calculation yields (as $k \ge 3$):
\begin{equation}\label{eq:2}
n = \frac {n^{k}}{(2b)^{k-1}} \left(\frac {2b}{n}\right)^{k-1} \le
\frac 14\frac {n^{k+1}}{(2b)^{k-1}}.
\end{equation}

By Claim~\ref{clm:finalstage} we have that after Avoider's $\ell$th
move there are at least $\ell - g - 1$ vertices of $A$-degree $k-1$
in $I_{\ell-1}$. Since every such vertex creates at least
$\left(\frac{3}{4} - o(1)\right)n$ unique threats, by using
\eqref{eq:1} and \eqref{eq:2} we get that the number of threats
after Avoider's $\ell$th move is at least
$$\left(\left(1-o(1)\right)\frac{27}{20} n\left(\frac
{n}{2b}\right)^{k-1} - 2\right) \left(\frac 34 - o(1)\right)n\ge
\frac {n^{k+1}}{(2b)^{k-1}} - \frac 32 n \ge \frac 58\frac
{n^{k+1}}{(2b)^{k-1}}.$$

As already mentioned, if the maximal degree in $A_\ell$ is less than
$k$, then Enforcer wins if and only if there are at least $r$
threats after Avoider's $\ell$th move.

By definition of $f^+$ and $b^+_{n,k}$, and since $b^+_{n,k} =
\omega(n)$ by Claim~\ref{clm:lowerBounds}, it is clear that
$b^+_{n,k} \le f^+_{\cK_{\cS_k}}$. However, in order to show that
$b^-_{n,k} \le f^-_{\cK_{\cS_k}}$ holds as well, we must show in
addition that Enforcer has a winning strategy if $b = O(n)$. Indeed,
if $b = o(n)$ Enforcer wins no matter how he plays since Avoider
will have $\omega(n)$ edges in his final graph, so assume $b =
\Theta(n)$. Enforcer does the following: before the game starts he
chooses an arbitrary set $U \subseteq V$ of size
$|U|=(2b)^{\frac{k}{k+1}} < n$, and in each step he claims some
arbitrary free edge with at least one endpoint outside $U$ until he
can no longer do so, i.e.\ until all free edges lie completely inside
$U$. Then he pretends to start a new game on $n'
=(2b)^{\frac{k}{k+1}}$ vertices with bias $b = \frac 12
n'^{\frac{k+1}{k}}$ according to the strategy for the case
$b=\omega(n)$. This is not exactly a new game because there may be
some edges inside $U$ already claimed by Avoider, and the ``new"
game may start during Enforcer's move. However, since Avoider can
claim only a constant number of edges incident to each vertex, and
since Enforcer makes at most $b$ additional steps before Avoider's
first move in the new game, these factors have no significant
effect. They only affect the case $j=1$ in the proof of
Claim~\ref{clm:avgDegreeStrict}, and it is easy to see that the
analysis there is still valid. The number of free edges before the
last round (i.e.\ $r(n',b)$) is also affected, but since $b \le
b^-_{n',k}$ Enforcer wins regardless of the exact value of $r$.
\end{proof}

\section{An application of Fact 2.1}\label{sec:apl}
We mentioned in the introduction that Bednarska-Bzd\c{e}ga in \cite{MBed} obtained
bounds on the different threshold biases for general $H$-games. For
every fixed graph $H$ with at least two edges, let
\[
m(H) = \max_{F\subseteq H: v(F)\ge 1} \hskip 0.1truecm
\frac{e(F)}{v(F)}; \hskip 0.7truecm m'(H) = \max_{F\subseteq H:
v(F)\ge 1}\hskip 0.1truecm\frac{e(F)-1}{v(F)};
\]
\[
m_2(H) = \max_{F\subseteq H: e(F)\ge 2}\hskip 0.1truecm \frac{e(F) -
1}{v(F)-2},
\]
where $v(F)$ and $e(F)$ denote the number of vertices and number of
edges of the subgraph $F$, respectively. Bednarska-Bzd\c{e}ga proved
the following (\cite{MBed}, Theorems~1.9 and~1.10):

\begin{enumerate}[$(i)$]
\item $f^{mon}_{\cK_H} = O(n^{1/m'(H)})$ and $f^+_{\cK_H} =
O(n^{1/m'(H)})$;
\item $f^{mon}_{\cK_H} = \Omega(n^{1/m_2(H)}/\ln
n)$ and $f^-_{\cK_H} = \Omega(n^{1/m_2(H)}/\ln n)$;
\item $f^-_{\cK_H} = O(n^{1/m(H)}\ln n)$ always holds, and
$f^-_{\cK_H} = O(n^{1/m(H)})$ holds for infinitely many values of
$n$.
\end{enumerate}

In the proof of $(iii)$ she uses her general criterion for Avoider's
win and shows that for $f = cn^{1/m(H)}$ (for some constant $c \ge
1$), if $r(n,f') > f$, then Avoider wins the strict $(1:f')$
$H$-game played on the edges of $K_n$. She uses two number theoretical
facts to show that there always exists such an $f'$ satisfying $f'
\le 4m(H)f\ln f$, and that for infinitely many values of $n$ there
exists such an $f'$ satisfying $f' \le 2f$. She only considers the
case $m(H) > 1/2$ (the other case is trivial), and so by applying
Fact~\ref{fact:largeRemainder} we obtain the following corollary.

\begin{corollary}\label{cor:improved}
If $m(H) \le 1$ then $f^-_{\cK_H} = O(n^{1/m(H)})$ .
\end{corollary}

That is, we obtain the better bound for all $n$ in this case. Note
that $m(H) \le 1$ if and only if $H$ is a graph in which every
connected component contains at most one cycle.

Our results show that these bounds are far from being tight for the
star game, except for the ``improved" upper bound on $f^-_{\cK_H}$
given in Corollary~\ref{cor:improved}, where we got exactly the same
result. We proved the bound for $f^-_{\cK_{\cS_k}}$ in our paper
explicitly anyway, since it is straightforward to obtain by using
Fact~\ref{fact:largeRemainder} and our other arguments. The gaps in
$(i)$ and $(ii)$ are not very surprising, as these bounds are valid
for \emph{every} fixed graph $H$. However, at least the upper bound
on $f^{mon}_{\cK_H}$ cannot be improved in general, since it is
tight for the case $H=K_3$. In addition, the constant exponent
$1/m(H)$ in both bounds of $(iii)$, as well as in
Corollary~\ref{cor:improved}, cannot be improved in general, because
for $H = P_3$ we have $f^-_{\cK_{P_3}}=\Omega(n^{1/m(H)})$, as
observed by Bednarska-Bzd\c{e}ga herself in~\cite{MBed}. In this
paper we provided an infinite family of graphs for which this bound
is tight.

\section{Concluding remarks and open problems}\label{sec:rmrks}

In Section~\ref{sec:enfStar} we propose a very natural strategy for
Enforcer in the $k$-star game to enforce the appearance of a vertex
of large degree in Avoider's graph. In~\cite{GS}, Gebauer and
Szab\'{o} use a very similar approach; they study the change of the
average degree in Breaker's graph over some subset of vertices
during the game and show that it cannot get too large, and thus
Maker's graph has a large minimum degree. So, enforcing a large
average degree (and thus the maximal degree) in Avoider's graph over
a subset of vertices complements in a way the method of Gebauer and
Szab\'{o}.

In this paper we show that for every sufficiently large $n$ and
every $k \geq 3$, the threshold biases $f^-_{\cK_{\cS_k}}$ and
$f^+_{\cK_{\cS_k}}$ are not of the same order, thus supporting the conjecture of Hefetz et al.~from \cite{HKSS10}. In that paper they also
showed that $f^{mon}_{\cK_{H}}$ and $f^-_{\cK_{H^-}}$ are of the
same order for $H = K_3$; we showed the same for $H = \cS_k$.
Observe that $H^- = P_3$ for both $H = K_3$ and $H = \cS_3$, and so
$f^{mon}_{\cK_{H}}$ is of the same order in both cases. It would be
interesting to further investigate the relation between the monotone
$H$-games and the strict $H^-$-games and to determine whether indeed
there is a connection between the two. Note that for a general graph
$H$, the graph $H^-$ is not uniquely determined (unlike the $K_3$
and $\cS_3$ cases) and so for different choices of $H^-$ there are
different outcomes. Therefore, choosing the ``correct" $H^-$ must
also be considered.

In Theorem~\ref{thm:mainTheorem} we provided the bounds for $f^{mon}_{\cK_{\cS_k}}$, $f^{+}_{\cK_{\cS_k}}$ and $f^{-}_{\cK_{\cS_k}}$ that are tight up to a constant factor. We could actually get some better (tighter) bounds -- for
example, by refining Avoider's strategy we could show that the constant in the upper
bound in all three cases is $1 + \varepsilon$, for any
$\varepsilon > 0$, rather than 2 -- but since we could not close the
gap completely we decided to provide slightly weaker bounds with
simpler proofs. It would be nice to determine the exact values of $C_1$, $C_2$ and $C_3$ for which
$f^{mon}_{\cK_{\cS_k}} = C_1 n^{\frac{k}{k-1}}$, $f^+_{\cK_{\cS_k}}
= C_2 n^{\frac{k}{k-1}}$ and $f^-_{\cK_{\cS_k}} = C_3
n^{\frac{k+1}{k}}$. In addition, note that our results for the
$k$-star game only hold for a constant $k$. 

Let us comment on the case when $k=k(n)$ tends to infinity along with $n$. Clearly, as long as the bias $b$ satisfies $b\le (1-\varepsilon)\frac{n}{k}$, Avoider is doomed as at the end of the game even the average degree of his graph will be larger than $k$. Thus $(1-o(1))\frac{n}{k} \le f^-_{\cK_{\cS_k}},f^+_{\cK_{\cS_k}}, f^{mon}_{\cK_{\cS_k}}$ holds. On the other hand, Avoider could win when $b\ge (1+\varepsilon)\frac{n}{k}$ provided he could keep all degrees asyptotically the same. This kind of \textit{discrepancy} games were studied first by Erd\H os, and the following general result of Beck tells us the order of magnitude of the threshold biases when $k=\omega(\log n)$. He considered the game of Balancer (playing with bias $p$) and Unbalancer (playing with bias $q$) in which they claim elements of a board $X$. 

\begin{theorem}[Theorem 17.5 in \cite{Be}]
\label{thm:balancer} Let $\cF$ be an arbitrary $N$-uniform hypergraph. Balancer and Unbalancer play the $(p:q)$ game: they alternate, Balancer takes $p$ new points and
Unbalancer takes $q$ new points per move. Then Balancer, as the first player, can
force that, at the end of the play, for every $A \in \cF$, his part in $A$ is strictly between $\frac{p+\varepsilon}{p+q} N$ and $\frac{p-\varepsilon}{p+q} N$, where
\[
\varepsilon= \left( 1+ O\left(pq\sqrt{\frac{\log |\cF|}{(p+q)N}}\right)\right)2pq\sqrt{\frac{\log |\cF|}{(p+q)N}}.
\]
\end{theorem}
 
Avoider can use Balancer's strategy in the above game with $p=1$, $q\ge (1+\varepsilon)\frac{n}{k}$, $N=n-1$ and $|\cF|=n$ with $\cF$ consisting of the edge sets of the stars of $K_n$. An easy calculation shows that if $k=\omega(\log n)$, then $\varepsilon$ can be chosen to satisfy $\varepsilon=o(1)$. Thus we obtain $f^-_{\cK_{\cS_k}},f^+_{\cK_{\cS_k}}=(1+o(1))\frac{n}{k}$ for these values of $k$. The case $\omega(1)=k=O(\log n)$ remains open.

\textbf{Acknowledgements:} The research was initiated at the 4th
Eml\'ekt\'abla Workshop held in Tihany, August 6-9, 2012. The
authors would like to thank Asaf Ferber and Milo\v{s} Stojakovi\'c
for their useful comments and suggestions. Also, the authors would
like to express their gratitude to the referees whose thorough remarks
greatly improved the paper. Furthermore, Fact 2.2 was pointed out by one of the referees along with the remark about the case $k \rightarrow \infty$ 
in the last section.


\begin{thebibliography}{99}

\bibitem{Be}
J. Beck, \textbf{Combinatorial Games: Tic-Tac-Toe Theory}, Cambridge University Press, 2008

\bibitem{MBed} M. Bednarska-Bzd\c{e}ga, Avoider-Forcer games on hypergraphs with small rank, \textit{Electronic Journal of Combinatorics}, 21(1) (2014), P1.2.

\bibitem{Boll}
B. Bollob\'{a}s, \textbf{Random Graphs}, Cambridge University Press,
2001.

\bibitem{CE}
V. Chv\'atal and P. Erd\H{o}s, Biased positional games, {\em Annals
of Discrete Mathematics} 2 (1978), 221--228.

\bibitem{GS}
H. Gebauer and T. Szab\'o, Asymptotic random graph intuition for the biased connectivity game, \emph{Random Structures and Algorithms},
35 (2009), 431--443.

\bibitem{HKSS10} D. Hefetz, M. Krivelevich, M. Stojakovi\'c and T. Szab\'o,
Avoider-Enforcer: The rules of the game, \textit{Journal of Combinatorial
Theory Series A} 117 (2010), 152--163.

\bibitem{HKS07} D. Hefetz, M. Krivelevich and T. Szab\'o,
 Avoider-Enforcer games, \textit{Journal of Combinatorial Theory Series
A} 114 (2007), 840--853.

\bibitem{JLR} S. Janson, T.\ \L uczak, A.\ Ruci\'{n}ski, \textbf{Random Graphs}, John Wiley \& Sons, Inc., 2000.

\bibitem{KS}
M.\ Krivelevich and T.\ Szab\'o, Biased positional games and small hypergraphs with large covers, \textit{Electronic Journal of Combinatorics}, 15 (2008), R70.

\bibitem{West} D. B. West, {\bf Introduction to Graph Theory}, Prentice Hall,
2001.
\end{thebibliography}
\end{document}